\documentclass[12pt]{amsart}
\usepackage{amssymb, mathtools,multirow}
\usepackage{multicol}
\usepackage{colordvi}
\usepackage{graphicx}
\usepackage{epsfig, epsf}
\numberwithin{equation}{section}
\theoremstyle{plain}
\newtheorem{Th}{Theorem}[section]
\newtheorem{Cor}[Th]{Corollary}
\newtheorem{Prop}[Th]{Proposition}

\numberwithin{equation}{section} \theoremstyle{definition}
\newtheorem{Rem}[Th]{Remark}

\newtheorem{Ex}[Th]{Example}
 
\newtheorem{Alg}[Th]{Algorithm}
\newtheorem{Def}[Th]{Definition}

\copyrightinfo{}{}
\newcounter{FNC}[page]
\def\fauxfootnote#1{{\addtocounter{FNC}{2}\Magenta{$^\fnsymbol{FNC}$}%
     \let\thefootnote\relax\footnotetext{\Magenta{$^\fnsymbol{FNC}$#1}}}}

\newcommand{\Z}{\mathbb Z}
\newcommand{\R}{\mathbb R}

\newcommand{\conv}{\operatorname{conv}}
\newcommand{\w}{\operatorname{w}}

\newcommand{\ls}{\operatorname{ls_\Delta}}

\newcommand{\lss}{{\rm ls}_\square}

\title[Lattice Size and Generalized Basis Reduction in Dimension 3]{Lattice Size and Generalized Basis Reduction in Dimension 3}

\author{Anthony Harrison}
\address{Department of Mathematics\\
        Kent State University\\
        1300 Lefton Esplanade, Kent, OH 44242\\
        USA}
\email{aharri60@kent.edu}
\thanks{}

\author{Jenya Soprunova}
\address{Department of Mathematics\\
        Kent State University\\
       1300 Lefton Esplanade, Kent, OH 44242\\
        USA}
\email{soprunova@math.kent.edu}
\urladdr{http://www.math.kent.edu/~soprunova/}

\subjclass[2010]{52B20, 11H06, 52C05, 52C07}
\keywords{Lattice size, successive minima, generalized basis reduction.}

\begin{document}

\begin{abstract}
The lattice size  of a lattice polytope $P$ was defined and studied by Schicho, and Castryck and Cools.
They provided an ``onion skins'' algorithm for computing the lattice size of a lattice polygon $P$ in $\R^2$ based on passing successively to the convex hull of the interior lattice points of $P$.

We explain the connection of the lattice size to the successive minima of $K=\left(P+(-P)\right)^\ast$ and to the lattice reduction with respect to the general norm that corresponds to $K$. 
It follows that the generalized Gauss algorithm of Kaib and Schnorr (which is faster than the ``onion skins'' algorithm) computes the lattice size of any convex body in $\R^2$.
 
 We extend the work of Kaib and Schnorr to dimension 3, providing a fast algorithm for lattice reduction with respect to the general norm defined by a convex origin-symmetric body $K\subset\R^3$.
  We also explain how to recover the successive minima of $K$ and the lattice size of $P$ from the obtained reduced basis and therefore provide a fast algorithm for computing  the lattice size of any convex body $P\subset\R^3$.

 \end{abstract}
\maketitle

\section*{Introduction}
The lattice size was defined and studied by  Castryck and Cools in~\cite{CasCools} in the context of simplification of a defining equation of an algebraic curve. Their work builds on the work of Schicho~\cite{Schicho}  who developed an ``onion skins'' procedure for unimodularly mapping a lattice polygon inside a small multiple of the standard simplex.  
 The {\it lattice size} $\operatorname{ls}_X(P)$ of a lattice polytope $P\subset\R^d$  with respect to a set $X\subset \R^d$ of positive Jordan measure is the smallest integer $l\geq 0$ such that $T(P)$ is contained in the $l$-dilate  of $X$ for some unimodular transformation $T$.  See Example~\ref{e:lsdef} for an illustration of this definition.
 
When $X=[0,1]\times\R^{d-1}$,  the lattice size of a lattice polytope $P$ with respect to $X$ is the lattice width $\w(P)$, an important invariant in lattice geometry and its applications, see the introduction to~\cite{CFB}. Of particular interest are also the cases when $X$ is either the standard $d$-dimensional simplex $\Delta$  or the unit cube $\square=[0,1]^d$, as explained at the end of this section.

The idea of the ``onion skins'' algorithm of~\cite{CasCools, Schicho} for computing the lattice size of a polygon $P$ is based on the observation that when one passes from $P$ to the convex hull of its interior lattice points, its lattice size $\ls(P)$ drops by 3 (and $\lss (P)$ by 2) unless $P$ belongs to a list of exceptional cases, which Castryck and Cools describe in~\cite{CasCools}.  One can then compute the lattice size by successively peeling off ``onion skins'' of the polygon. 

Let $P\subset\R^ d$ be a convex body and consider $K=(P+(-P))^\ast$, the polar dual of the Minkowski sum of $P$ with $-P$, its reflection in the origin.  In Section 1 we show  that the lattice size $\lss(P)$ with respect to the unit cube can be recognized as a version of the Minkowski successive minimum of $K$ (see Definitions~\ref{D:Successive_min1} and~\ref{D:Successive_min2} ). Computing the successive minima of $K$ is related to finding a lattice basis of vectors which are short with respect to the general norm defined by $K$. The problem of finding a lattice basis of short vectors, a so-called reduced basis,  is a hard and important problem that has seen a lot of study, in particular due to its applications in cryptology (see, for example, ~\cite{Barvinok} and Chapter IV of~\cite{Gal}). 
The LLL algorithm~\cite{LLL} of A. Lenstra, H. Lenstra, and Lov\'asz from 1982 was an important breakthrough in the theory of lattice reduction in the case of inner-product norm. Since then, other algorithms have been developed, with most of them dealing with the case of the inner product norm.

In~\cite{LovScarf} Lov\'asz and Scarf generalized the LLL reduction algorithm to the case of the general norm defined by an origin-symmetric convex body $K$. They showed that the algorithm takes polynomial time when the dimension is fixed, and that the norms of the output basis vectors provide a bound for  the successive minima of $K$. In the case of dimension 2, the algorithm (called the generalized Gauss algorithm) was analyzed by Kaib and Schnorr in~\cite{KaibSch}. They showed that  the number of iterations of the algorithm is logarithmic in $\max\{\Vert h^1\Vert,\Vert h^2\Vert\}/\lambda_2$, where $(h^1,h^2)$ is the input lattice basis and $\lambda_2$ is the second successive minimum of $K$.
They also showed that the norms of the output basis are equal to the successive minima of $K$.  Hence, in particular, the generalized Gauss algorithm can be used to compute the lattice size $\lss(P)$ of a convex body $P\subset\R^2$. This algorithm is faster than  the ``onion skins'' algorithm of~\cite{CasCools}, as we explain after Theorem~\ref{T:main2D} below. The algorithm also works for any convex body $P\subset\R^2$, while the ``onion skins'' algorithm can only be used in the case when $P$ is a lattice polygon.

Although we were initially motivated by the problem of computing the lattice size of a lattice polytope, the main result of this paper is the analysis of the generalized lattice basis reduction algorithm in $\R^3$, which extends the results  of~\cite{KaibSch} to dimension~3. This provides  a fast algorithm for computing both, the successive minima  of an origin-symmetric convex body $K\subset\R^3$, and  the lattice size of an arbitrary convex body $P\subset\R^3$. In~\cite{KaibSch} the authors used generalized continuants for their analysis of the algorithm in $\R^2$. This approach does not seem to be feasible to extend to dimension 3. Instead, in Section 2, we first modify the algorithm in dimension 2 to make the analysis easier and then, in Section 3, extend our work to dimension 3. (See the comment before Algorithm~\ref{A:alg2} explaining the nature of the modification.) We show that the number of iterations in dimension 3 is  logarithmic in the sum of norms of the input basis vectors, divided by the sum of the successive minima of $K$. We also show that the algorithm outputs a reduced basis with one of the norms of the basis vectors equal to $\lss(P)$, and explain how to recover the successive minima of $K$. Recently, we showed in~\cite{HST}  that the generalized lattice basis reduction can also be used to compute the lattice size $\ls(P)$ of a plane convex body $P$ with respect to the standard simplex $\Delta$, and provide an example demonstrating that this cannot be extended to dimension 3.

In the classical case of the inner product norm, the basis reduction also simplifies significantly in the case of low dimension, as explained in~\cite{Nguen}. The case of dimension~3 is considered in~\cite{Semaev, Vallee}.  

In~\cite{CasCools, Schicho} the motivation for studying the lattice size comes from algebraic\linebreak geometry.  Let $k$ be a field and consider a hypersurface in $\left(k\setminus\{0\}\right)^d$ defined by  a Laurent polynomial $f\in k[x_1^{\pm1},\dots, x_d^{\pm 1}]$. The {\it Newton polytope} $P=P(f)$ of $f$ is  the convex hull of the exponent vectors that appear in $f$. The total degree of $f$ can then be interpreted as the smallest $l$ such that $P$ is contained in the $l$ dilate $l\Delta$ of the standard simplex $\Delta$ after a shift by a lattice vector. Let $A=(a_{ij})$ be a unimodular matrix, that is, all entries in $A$ are integers and $\det A=\pm 1$. Then $A$ defines a monomial change of variables  
$x_i=u_1^{a_{1i}}u_2^{a_{2i}}\dots u_d^{a_{di}}$  after which $f$ turns into a Laurent polynomial with Newton polytope $A(P)$.  Hence $\ls(P)$ is the smallest total degree of a Laurent polynomial obtained from $f$ by making such monomial changes of variables. 
 
 This monomial change of variables is a particular case of a birational map of the hypersurface defined by $f=0$ in the algebraic torus  $(k\setminus \{0\})^d$. Hence $\ls(P)$ provides an upper bound for the lowest total degree of the hypersurface defined by $f=0$ in $(k\setminus \{0\})^d$ under birational equivalence. Note that when finding $\lss(P)$, the lattice size of $P$ with respect to the cube,  we are minimizing  the largest component of all the multi-degrees of the monomials that appear in $f$,  over monomial changes of variables.
 
Lattice size is also useful when dealing with questions that arise when studying lattice polytopes. In \cite{BrownKasp} the lattice size of a lattice polygon with respect to the unit square is used to classify small lattice polygons and corresponding toric codes.  This notion also appears implicitly in \cite{Arnold}, \cite{Barany_Pach}, and \cite{LagZieg}.

\smallskip

\noindent {\it Data Availability Statement.} Data sharing is not applicable to this article as no datasets were generated or analyzed during the current study.


\section{Relation to successive minima}

Let $P$ be a convex body in $\R^d$ with non-zero volume.

\begin{Def}
A vector $(a_1,\dots,a_d)\in\Z^d$ is primitive if $\gcd(a_1,\dots, a_d)=1.$
\end{Def}

\begin{Def}
Given $h\in\Z^d$, the {\it lattice width of $P$ in the direction of $h$} is defined by
$$\w_h(P)=\max_{x\in P}(x\cdot h) -  \min_{x\in P} (x\cdot h),
$$
where $x\cdot h$ denotes the standard dot product.  The {\it lattice width} $\w(P)$  is the minimum of $\w_h(P)$ over all non-zero primitive integer directions $h\in\Z^d$.  
\end{Def}

Recall that a square integer matrix $A$ is called unimodular if $\det A=\pm 1$. For such a matrix $A$, a transformation $T\colon \R^d\to\R^d$, defined by $T(x)=Ax+v$, where $v\in\Z^d$, is called an {\it affine unimodular transformation}. Such transformations preserve the integer lattice $\Z^d\subset\R^d$.

\begin{Def}
The {\it lattice size} $\operatorname{ls}_X(P)$ of a non-empty convex body $P\subset\R^d$  with respect to a set $X$ of positive Jordan measure is the smallest $l\geq 0$ such that $T(P)$ is contained in the $l$-dilate  of $X$ for some affine unimodular transformation $T$ of the form $T(x)=Ax+v$. We say that such $A$ {\it computes} the lattice size.
\end{Def}

Let $\square=[0,1]^d\subset \R^d$ denote the unit cube in $\R^d$. In this paper we will be studying  $\lss(P)$, the lattice size of $P$ with respect to the unit cube.

%
%
\begin{Ex}\label{e:lsdef} Let $P\subset\R^2$ be the triangle with vertices $(0,0)$, $(1,0)$, and $(2,3)$. 
\vspace{.2cm}
\hspace{3cm} \includegraphics[scale=1]{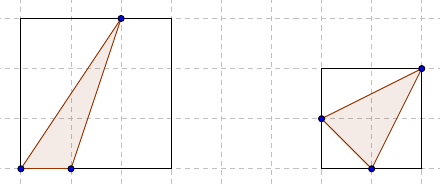} 

\noindent To get from the first  to the second polygon in this  diagram we apply the map $$\left[\begin{matrix}x\\y\end{matrix}\right]\to\left[\begin{matrix}1&0\\-1&1\end{matrix}\right]\left[\begin{matrix}x\\y\end{matrix}\right]+\left[\begin{matrix}0\\1\end{matrix}\right].$$
Since $P$ contains an interior lattice point, it is impossible to unimodularly map $P$ inside the unit square and hence $\lss(P)=2$. Also, $A=\left[\begin{matrix}1&0\\-1&1\end{matrix}\right]$ computes the lattice size of $P$.
\end{Ex}
\begin{Def} For $k=1,\dots,d$ we define $\operatorname{ls}_k(P)$ to be the lattice size of $P$ with respect to $X=[0,1]^k\times\R^{d-k}$.
Note that $\operatorname{ls}_d(P)=\lss(P)$ and $\operatorname{ls}_1(P)=\w(P)$.
\end{Def}

Given a convex origin-symmetric body $K\subset\R^d$, the norm associated with $K$ is defined by  
$$\Vert h\Vert _K\colon=\inf\{\lambda>0\mid h \in \lambda K\},\ {\rm for}\ h\in\R^d.$$
This indeed defines a norm on $\R^d$ and any norm on $\R^d$ is of this form for some convex origin-symmetric $K$. In this case $K=\{h\in\R^d\mid \Vert h\Vert_K\leq 1\}$, that is, $K$ is the unit ball with respect to the norm.

We next define the successive minima of $K$. They were first defined and studied by Minkowski and since then have played an important role, in particular, in geometry of numbers and computational number theory~\cite{Gruber, Siegel}.
Let $K$ be a convex origin-symmetric body and let $\Vert \cdot\Vert _K$ be the norm associated with $K$.

\begin{Def}\label{D:Successive_min1} 
We say that $\lambda_1,\dots,\lambda_d$ are the {\it (Minkowski) successive minima of $K$} if  $\lambda_k$   is the smallest positive real number such that $\lambda_kK$ contains $k$ linearly independent lattice vectors, where $k=1,\dots, d$. 
That is, $\lambda_k=\Vert h^k\Vert_K $, where  $h^1,\dots, h^k\in\Z^d$  and  $h^{k}$  is the shortest lattice vector linearly independent of $h^1, h^2,\dots, h^{k-1}$.
 \end{Def}

\begin{Ex} Here is a simple example illustrating this definition. Consider $K=\conv\{(\pm2,0), (0,\pm 1)\}\subset\R^2$. Then $\lambda_1(K)=\frac{1}{2}$ and  $\lambda_2(K)=1$.

\vspace{.2cm}
\begin{center}
 \includegraphics[scale=1]{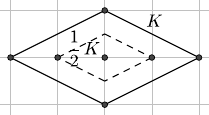} 
 \end{center}
\end{Ex}

The lattice vectors that correspond to the successive minima do not necessarily form a basis of the lattice, as will be demonstrated in Example~\ref{E:ex_mu}.  
Larger invariants that correspond to lattice vectors that do form a basis of the lattice were considered in~\cite{Siegel}. 
 We next define this version of the successive minima.

\begin{Def}\label{D:Successive_min2}
We define $\mu_1,\dots\mu_d$, so that each $\mu_k$ is the smallest positive real number such that $\mu_kK$ contains $k$ linearly independent lattice vectors that can be extended to a basis of the 
integer lattice $\Z^d$.  That is, $\mu_k=\Vert h^k\Vert_K $, where  each $h^{k}$  is the shortest lattice vector 
such that $(h^1,\dots, h^k)$ can be extended to a basis of $\Z^d$.
\end{Def}

Note that $\mu_i$ are the norms of the basis vectors in a Minkowski reduced lattice basis introduced by H. Minkowski, see Definition~\ref{D:Minkowski_reduced} below.

Clearly, $\lambda_k\leq \mu_k$. Furthermore, $\lambda_1=\mu_1$ and, due to Pick's formula, we also have $\lambda_2=\mu_2$.

\begin{Ex}\label{E:ex_mu} Consider $K=\conv\{(\pm 1,0,0), (0,\pm 1,0), \pm(1,1,2)\}\subset\R^3$. Note that the lattice points of $K$ are its vertices listed in the definition of $K$, and the origin.
We have $\lambda_1(K)=\lambda_2(K)=\lambda_3(K)=1$. Note that $K$ does not contain a basis of $\Z^3$, but its multiple $\frac{3}{2}K$  has six additional lattice points, 
$\pm(1,0,1), \pm(0,1,1)$, and $\pm(1,1,1)$, all of them on the boundary. We conclude that $\mu_3(K)=\frac{3}{2}$.
\end{Ex}

We next explain an important connection between the lattice size of $P$ and the successive minima of  $K=(P+(-P))^\ast$, the polar dual of the Minkowski sum of $P$ with its refection in the origin $-P$. 


\begin{Th}\label{T:translate} 
\begin{itemize}
\item[(a)] For  $h\in\Z^d$ we have  $\Vert h\Vert_K=\w_h(P)$.
\item[(b)] We have $\operatorname{ls}_k(P)=\mu_k(K)$ for $k=1\dots d$. In particular, $\lss(P)=\mu_d(K)$.
\end{itemize}
\end{Th}

\begin{proof} To prove part (a) we observe that
 \begin{eqnarray*}
\Vert h\Vert_K&=&\inf\{\lambda>0\mid h/\lambda \in K\}\\
&=&\inf\{\lambda>0\mid \frac{h}{\lambda}\cdot x\leq 1 {\rm\ for\ all\ } x\in K^\ast\}=\max\limits_{x\in K^{\ast}}  h\cdot x\\
&=&\frac{1}{2}\w_h(K^\ast)=\w_h(P).
\end{eqnarray*}


For part (b) recall that $\operatorname{ls}_k(P)$ is the lattice size of $P$ with respect to $X=[0,1]^k\times\R^{d-k}$. Therefore, for some unimodular matrix $A$ and a lattice vector $v$, we have 
$$A(P+v)\subset \operatorname{ls}_k(P)[0,1]^k\times\R^{d-k}.$$ Hence $A(P+(-P))\subset \operatorname{ls}_k(P)[-1,1]^k\times\R^{d-k}$. Taking the polar dual of both sides we get
$$(A^{T})^{-1}(P+(-P))^\ast\supset \frac{1}{\operatorname{ls}_k(P)}\Diamond_k\times\{0\}.$$
Here $\Diamond_k$ is the cross-polytope, that is,
$\Diamond_k=\conv\{\pm e^i\mid i=1,\dots, k\}$, where $e^i$ are the standard basis vectors. This implies that $\operatorname{ls}_k(P)K\supset A^{T}(\Diamond_k)$.
Since $K$ is symmetric with respect to the origin, $\operatorname{ls}_k(P)$ is the smallest positive number such that $\operatorname{ls}_k(P)K$ contains $k$ linearly independent vectors that are unimodularly equivalent to 
$(e^1,\dots, e^k)$.
We conclude that $\operatorname{ls}_k(P)=\mu_k(K)$ and,  in particular, $\lss(P)=\mu_d(K)$.
\end{proof}

We will observe in Proposition~\ref{P:Acomputes} that matrix $A$ whose rows are  vectors of a Minkowski reduced basis (see Definition~\ref{D:Minkowski_reduced}) computes  the lattice size.

We have recognized the lattice size as a successive minimum and therefore we have reduced our initial question of computing the lattice size to finding the successive minima of an origin-symmetric convex body $K$. 

 In~\cite{LovScarf} Lov\'asz and Scarf  provided a basis reduction algorithm that bounds the successive minima. In dimension 2 the algorithm was analyzed by Kaib and Schnorr in~\cite{KaibSch}. They estimated the complexity of the algorithm and showed that the successive minima are equal to the norms of the output vectors. Therefore, this algorithm also computes the lattice size of an arbitrary convex body in the plane.

In Chapter 3 we extend the work of Kaib and Schnorr to dimension three and therefore provide a fast basis reduction algorithm  that computes the successive minima of an origin-symmetric body $K\subset\R^3$ as well
the lattice size of a convex body $P\subset\R^3$.

\section{Basis Reduction and  Successive Minima in Dimension 2.}

In this section we deal with computing the successive minima of a convex origin-symmetric body $K\subset\R^2$. This is accomplished using the generalized Gauss algorithm for basis reduction, explained and analyzed in~\cite{KaibSch}.  In order to extend the algorithm to dimension 3, we will first adjust the 2-dimensional algorithm of~\cite{KaibSch} to make its analysis easier.
Our methods are different than the ones used in \cite{KaibSch}.

Note that if we start with a convex body $P\subset\R^2$ we then work with $K=\left(P+(-P)\right)^\ast\subset\R^2$.  If we start with $K$, an arbitrary convex origin-symmetric body in $\R^2$, we can pick $P=\frac{1}{2}K^{\ast}$ and then  $K=\left(P+(-P)\right)^\ast$. Therefore, even if starting with $K$, we can always refer to a corresponding $P$. 

Let $K$ be a convex origin-symmetric body in $\R^2$  and consider the lattice $\Z^2\subset \R^2$. Let $\Vert \cdot\Vert $ be the norm on $\Z^2$ associated with $K$, that is,
for $h\in\Z^2$ let  $\Vert h\Vert=\Vert h\Vert_K=\w_{h}(P)$. 

\begin{Def}\label{D:red2D} We say that a basis $(h^1,h^2)$ of the integer lattice $\Z^2$ is {\it reduced} (with respect to $K$), if 
\begin{itemize}
\item[(1)] $\Vert h^1\Vert \leq \Vert h^2\Vert$ and
\item[(2)]  $\Vert h^1\pm h^2\Vert \geq \Vert h^2\Vert$.
\end{itemize}
\end{Def}

\begin{Ex}
Let $P$ be the first triangle from Example~\ref{e:lsdef} and $K=\left(P+(-P)\right)^\ast$.  Let $h^1=(1,0)$ and $h^2=(-1,1)$. We have 
$\Vert h^1\Vert =\w_{(1,0)}(P)=2$, $\Vert h^2\Vert =\w_{(-1,1)}(P)=2$, $\Vert h^1+h^2\Vert =\w_{(0,1)}(P)=3$, $\Vert h^1-h^2\Vert =\w_{(2,1)}(P)=7$, and hence $(h^1,h^2)$ is a reduced basis of $\Z^2$.
\end{Ex}

Note that part $(2)$ in Definition~\ref{D:red2D} can be replaced  with a seemingly stronger condition $\Vert m h^1+ h^2\Vert \geq \Vert h^2\Vert$ for all $m\in\Z$ since for $m>0$ we have
$$\Vert m h^1\pm h^2\Vert + (m-1)\Vert  h^2\Vert\geq m\Vert h^1\pm h^2\Vert \geq m \Vert h^2\Vert,
$$
which implies that $\Vert m h^1\pm h^2\Vert\geq  \Vert h^2\Vert$.

Let $(h^1,h^2)$ be a reduced basis of $\Z^2$ with respect to the norm defined by $K$. We explain next that $h^1$ and $h^2$ have norms equal to the successive minima of $K$. Hence to compute the successive minima of $P\subset\R^2$  it is enough to find a basis which is reduced with respect to $K$. This theorem is proved in~\cite{KaibSch}, but we provide a simpler argument. 

\begin{Th}\label{T:2Dmain} Let $K$ be a convex origin-symmetric body in $\R^2$. If $(h^1,h^2)$ is a reduced basis of $\Z^2$ with respect to the norm defined by $K$, then $\lambda_1(K)=\Vert h^1\Vert$ and $\lss(P)=\lambda_2(K)=\Vert h^2\Vert$, where 
$\Vert \cdot\Vert $ is the norm defined by $K$.
\end{Th}

\begin{proof} It is enough to show that $\Vert ah^1+bh^2\Vert\geq \Vert h^2\Vert$ for all primitive $(a,b)$ except for $(1,0)$.
Switching, if necessary, the signs of $h^1$ and $h^2$, we can assume that $a, b>0$. If $a\geq b$ we get 
$$\Vert ah^1+bh^2\Vert + (a-b)\Vert h^2\Vert \geq a\Vert h^1+ h^2\Vert \geq a\Vert h^2\Vert ,
$$
which implies $\Vert ah^1+bh^2\Vert\geq b \Vert h^2\Vert\geq  \Vert h^2\Vert.$
Similarly, if $b\geq a$ we get 
$$\Vert ah^1+bh^2\Vert + (b-a)\Vert h^1\Vert \geq b\Vert h^1+ h^2\Vert \geq b\Vert h^2\Vert ,
$$
so $\Vert ah^1+bh^2\Vert\geq a \Vert h^2\Vert\geq  \Vert h^2\Vert.$
\end{proof}

%
%
%

Let $(h^1,h^2)$ be a basis of $\Z^2$ such that  $\Vert h^1\Vert\leq\Vert h^2\Vert$. Let $\min\limits_{m\in\Z}\Vert mh^1+h^2\Vert$ be attained at $f=mh^1+h^2$. 

\begin{Prop}\label{P:2Dreduced} If  $\Vert f\Vert\geq \Vert h^1\Vert$ then  the basis $(h^1,f)$ is reduced.
 \end{Prop}
 \begin{proof}  We have
 $$\Vert h^1\pm f\Vert = \Vert (m\pm 1)h^1+h^2\Vert\geq \Vert f\Vert.
$$
 \end{proof}
 
This proposition implies that the following algorithm (provided in~\cite{KaibSch}) computes a reduced basis of $\Z^2$.

\begin{Alg}\label{A:alg1} Let  $\Vert\cdot\Vert$ be the norm associated with an origin-symmetric convex body $K\subset\R^2$, and let $(h^1,h^2)$ be a basis  of $\Z^2$ that satisfies $\Vert h^1\Vert\leq\Vert h^2\Vert$. 
At each step, find $f$ such that the minimum of $\Vert mh^1+h^2\Vert$ over $m\in\Z$ is attained at $f=mh^1+h^2$. If $\Vert f\Vert< \Vert h^1\Vert$, pass to the basis $(f, h^1)$ and repeat the step. If $\Vert f\Vert\geq \Vert h^1\Vert$ then 
output a reduced basis $(h^1,f)$.
\end{Alg}

Let us show that this algorithm will terminate after finitely many steps.  For this we fix a number $l>0$ and show that there are only finitely many integer directions $h\in\Z^2$ such that 
$\Vert h\Vert_K \leq l$.  Consider a ball $B\subset\R^3$ such that $K\subset B$. Shifting $K$ we can assume that the center of $B$ is at the origin. Let $R$ be the radius of $B$ and let
 $\Vert h\Vert_2$ be the Euclidean norm of $h$.  Then if $\Vert h\Vert_2 \geq \frac{l}{R}$ we get
 $$\Vert h\Vert_{K}\geq \Vert h\Vert_{B}=R\cdot\Vert h\Vert_2 \geq l.
 $$

Since there are finitely many integer vectors with  $\Vert h\Vert_2 < \frac{l}{R}$, there are finitely many values of $\Vert h^2\Vert_K$ that the algorithm may traverse, so it will terminate after finitely many iterations.
(This argument is almost identical to the proof of Lemma 4.1  in~\cite{BrownKasp}.)

In~\cite{KaibSch}  the authors used generalized continuants to estimate the number of iterations the algorithm takes. Their approach does not seem to be feasible to extend to dimension 3. 
We next modify the algorithm to make the analysis easier so that we can later extend it to dimension 3. The main step of the new algorithm is similar to that of Algorithm~\ref{A:alg1}, but  we add a termination condition:
we now stop and output a reduced basis if  $\frac{2}{3}\Vert h^2\Vert\leq \Vert f\Vert$, which  allows us to estimate the number of iterations the algorithm takes.

 \begin{Alg}\label{A:alg2}
 Let  $\Vert\cdot\Vert$ be the norm associated with an origin-symmetric convex body $K\subset\R^2$. Start with a basis $(h^1,h^2)$ of $\Z^2$ that satisfies $\Vert h^1\Vert\leq\Vert h^2\Vert$. 
 At each step, find $f$ such that the minimum of $\Vert mh^1+h^2\Vert$ over $m\in\Z$ is attained at~$f=mh^1+h^2$. 
\begin{itemize}
\item[(i)] If $\Vert f\Vert <\Vert h^1\Vert$ and $\Vert f\Vert <\frac{2}{3}\Vert h^2\Vert$ then pass to the basis $(f,h^1)$ and repeat the step.
\item[(ii)] If $\Vert f\Vert\geq \Vert h^1\Vert$ then by Proposition~\ref{P:2Dreduced} we have found a reduced basis $(h^1,f)$.
\item[(iii)] If  $\frac{2}{3}\Vert h^2\Vert\leq \Vert f\Vert< \Vert h^1\Vert$ then output the two vectors of the smallest norm in the set $\{h^1, h^1\pm h^2, 2h^1\pm h^2\}$, which by Proposition~\ref{P:exceptional2D} below form a reduced basis.
\end{itemize}
\end{Alg}
 
Let us justify part (iii) of the algorithm.

\begin{Prop}\label{P:exceptional2D} Suppose that  $\frac{2}{3}\Vert h^2\Vert\leq \Vert f\Vert< \Vert h^1\Vert$ and consider the set $S=\{h^1, h^1\pm h^2, 2h^1\pm h^2\}$.
Then the two vectors of  the smallest norm in this set form a reduced basis. 
\end{Prop}

\begin{proof}
Let $(a,b)\in\Z^2$ be primitive with $b\geq 0$. We will show that $\Vert ah^1+bh^2\Vert  \geq \Vert h^2\Vert $ for all such $(a,b)$ except, possibly, for $(1,0), (\pm 1, 1)$, and  $(\pm 2, 1)$. 
Let  $b>0$ and let $kb$ be an integer multiple of $b$ closest to $a$. Then $|a-kb|\leq \frac{b-1}{2}$ for odd $b$ and $|a-kb|\leq \frac{b}{2}$ for even $b$.
For even $b\geq 4$ the latter inequality can be improved to $|a-kb|\leq \frac{b-2}{2}$ since $|a-kb|=\frac{b}{2}$ would  imply that both $a$ and $b$ are multiples of $b/2$, which is impossible as $(a,b)$ is primitive. 
Using the triangle inequality we get 
$$\Vert ah^1+bh^2\Vert +|a-kb|\Vert h^1\Vert\geq b \Vert kh^1+h^2\Vert\geq b \Vert f\Vert\geq  \frac{2 b}{3}\Vert h^2\Vert. 
$$
Hence for odd $b\geq 3$ we get $\Vert ah^1+bh^2\Vert\geq \left(\frac{2 b}{3}-\frac{b-1}{2}\right) \Vert h^2\Vert\geq \Vert h^2\Vert$. Similarly, for even $b\geq 4$ we get
$\Vert ah^1+bh^2\Vert\geq \left(\frac{2 b}{3}-\frac{b-2}{2}\right) \Vert h^2\Vert\geq \Vert h^2\Vert$.

For $b=2$ if $|a|\geq 5$ we have
\begin{equation*}
\Vert ah^1+2h^2\Vert \geq |a|\Vert h^1\Vert -2\Vert h^2\Vert\geq \left(\frac{2|a|}{3}-2\right)\Vert h^2\Vert\geq \Vert h^2\Vert.
\end{equation*}
If $|a|=4$ or $2$ then $(a,b)$ is not primitive. If $|a|=3$ then
$$\Vert 3h^1\pm 2 h^2\Vert+\Vert h^2\Vert\geq 3\Vert h^1\pm h^2\Vert\geq 3 \Vert f\Vert \geq 2\Vert h^2\Vert,$$ 
and hence $\Vert 3h^1\pm 2 h^2\Vert\geq \Vert h^2\Vert$. If $|a|=1$ we have $\Vert h^1\pm 2 h^2\Vert\geq 2\Vert h^2\Vert-\Vert h^1\Vert\geq \Vert h^2\Vert$.

Finally, if $b=1$ and $|a|\geq 3$ then $\Vert ah^1+h^2\Vert\geq   |a|\Vert h^1\Vert -\Vert h^2\Vert\geq  \Vert h^2\Vert$.
We have checked that $\Vert ah^1+bh^2\Vert  \geq \Vert h^2\Vert $ for all primitive $(a,b)$ with $b\geq 0$ except for $(\pm 2, 1)$, $(\pm 1, 1)$, and $(1,0)$.
Note that this also tells us that under our assumptions $f$ is one of $\{h^1\pm h^2, 2h^1\pm h^2\}$.

We have checked that to find a reduced basis it is enough to consider vectors in~$S$. Consider the two shortest vectors in $S$. The triangle formed by these two vectors cannot have any other points in $\Z^2$. Hence these two vectors form a basis of $\Z^2$, and this basis is reduced.
\end{proof}

When estimating the number of steps Algorithm~\ref{A:alg2} takes we will assume that we can use a norm oracle, that is, for a given $h\in\Z^2$ we can use the oracle output $\Vert h\Vert $.

\begin{Th}\label{T:main2D}
Algorithm~\ref{A:alg2} takes at most $\log_{\frac{6}{5}}\frac{\Vert h^1\Vert+\Vert h^2\Vert}{\lambda_1+\lambda_2}$ iterations and produces a reduced basis in 
$\mathcal{O}\left(\log \frac{\Vert h^2\Vert}{\lambda_1}\right)$
operations.
\end{Th}

\begin{proof}
At each iteration we pass from $(h^1, h^2)$  with $\Vert h^1\Vert\leq \Vert h^2\Vert$ to $(f,h^1)$ with  $\Vert f\Vert <\frac{2}{3}\Vert h^2\Vert$.
Hence 
$$\Vert h^1\Vert +\Vert f\Vert\leq \frac{5}{6}\left(\Vert h^1\Vert +\Vert h^2\Vert\right).
$$
Let $k$ be the number of iterations. Then the algorithm terminates in $k$ steps with a basis whose sum of norms is $\lambda_1+\lambda_2$, so we get $\lambda_1+\lambda_2\leq \left(\frac{5}{6}\right)^k \left(\Vert h^1\Vert +\Vert h^2\Vert\right)$
and the bound on the number of iterations follows.

Let's first consider a step which is not final, so at this step we pass from  $(h^1, h^2)$ to $(f,h^1)$, where $\Vert f\Vert<\Vert h^1\Vert$.
At each such step we need to minimize $\Vert mh^1+h^2\Vert$. We have
$$\Vert mh^1+h^2\Vert\geq |m|\Vert h^1\Vert-\Vert h^2\Vert\geq \Vert h^1\Vert,
 $$
 provided that $|m|\geq \frac{\Vert h^2\Vert}{\Vert h^1\Vert}+1$. Also,
 
$$\Vert mh^1+h^2\Vert\geq \Vert h^2\Vert - |m| \Vert h^1\Vert\geq \Vert h^1\Vert
$$
for $|m|\leq \frac{\Vert h^2\Vert}{\Vert h^1\Vert}-1$. Hence the only values of $m$ for which we may have \\ $\Vert mh^1+h^2\Vert<\Vert h^1\Vert$ are $\pm\left\lfloor \frac{\Vert h^2\Vert}{\Vert h^1\Vert}\right\rfloor$ and 
$\pm\left\lceil \frac{\Vert h^2\Vert}{\Vert h^1\Vert}\right\rceil$, where these four values reduce to two if $\Vert h^2\Vert$ is an integer multiple of $\Vert h^1\Vert$. If for at least one of these values of $m$ we have $\Vert mh^1+h^2\Vert \leq \Vert h^1\Vert$ then to find the minimum of $\Vert mh^1+h^2\Vert$ we only need to look at the values of the norm at these four values of $m$ and pick the smallest one, which takes $\mathcal{O}(1)$ steps.  If for all of these four values of $m$ we have $\Vert mh^1+h^2\Vert\geq \Vert h^1\Vert$, then this inequality holds true for all integer $m$, so we can conclude $\Vert f\Vert\geq \Vert h^1\Vert$ which tells as that we are at the final iteration. To minimize $\Vert mh^1+h^2\Vert$ at the final iteration, we observe that for $|m|\geq \frac{2\Vert h^2\Vert}{\Vert h^1\Vert}$ we have 
$$\Vert mh^1+h^2\Vert\geq |m|\Vert h^1\Vert-\Vert h^2\Vert\geq \Vert h^2\Vert.
 $$
Hence to find the minimum at the  last iteration  we need to compare the values of  $\Vert mh^1+h^2\Vert$ at $m\in\Z$ that satisfy $|m|<\frac{2\Vert h^2\Vert}{\Vert h^1\Vert}\leq \frac{2\Vert h^2\Vert}{\lambda_1}$.
Using ternary search, the last step can be completed in $\mathcal{O}\left(\log \frac{\Vert h^2\Vert}{\lambda_1}\right)$ operations.
Since the number of iterations is of order  $\log\frac{\Vert h^1\Vert+\Vert h^2\Vert}{\lambda_1+\lambda_2}\leq \log\frac{\Vert h^2\Vert}{\lambda_1}$, the complexity of the algorithm is $\mathcal{O}\left(\log \frac{\Vert h^2\Vert}{\lambda_1}\right)$.
\end{proof}

Note that in the case when $P$ is a lattice polygon, finding the lattice width in a  given direction requires $\mathcal{O}(v)$ operations, where $v$ is the number of vertices of $P$.  Since at most two vertices of $P$ may have the same $y$-coordinate, we have
$v\leq 2\lambda_1+2$ and hence Algorithm~\ref{A:alg2} computes the lattice size of $P$ in $\mathcal{O}\left(\lambda_1\log \frac{\Vert h^2\Vert}{\lambda_1}\right)$ operations.
There is no computation of the running time of the ``onion skins'' algorithm given in~\cite{CasCools}, but the implementation of the algorithm is provided in the accompanying Magma files. The algorithm lists the interior lattice points of $P$ and then iteratively applies this operation to the convex hull of the interior lattice points of $P$. In our notation, each such iteration requires $\mathcal {O}(\Vert h^1\Vert\Vert h^2\Vert)$ operations and there are   $\mathcal {O}(\lambda_1)$ of them. We conclude that Algorithm~\ref{A:alg2} outperforms the ``onion skins'' algorithm of~\cite{CasCools}.

Before we provide an example of how the algorithm works, let's establish how the change of basis matrix $A$ should be applied to  $K$ and $P$. Let the columns of $A$ be $v_1,\dots,v_d$ and let $e^1,\dots,e^d$ be 
the standard basis vectors. It follows directly from the definition that for $h\in\R^d$ we have
$$\Vert A h\Vert_K =\Vert h\Vert_{A^{-1}K}\ {\rm and} \w_{Ah}(P)=\w_h(A^{T}P).
$$
It follows that $\Vert e^i\Vert_{A^{-1}K}=\Vert Ae^i\Vert_K=\Vert v_i\Vert_K$ \ {\rm and} $\w_{e^i}(A^{T}P)=\w_{Ae^i}(P)=\w_{v_i}(P)$. Hence we should be applying $A^{-1}$ to $K$ and $A^{T}$ to $P$.

\begin{Ex} Let  $f(x)=4x^4-22x^3+49x^2y^2-52xy^3+22y^4$ and let  $K\subset\R^2$ be the level set defined by $f(x,y)\leq 1$. It is clear that $K$ is origin symmetric. We will explain why $K$ is convex after we find a reduced basis.
Let $(m,n)\in\Z^2$ and let $\lambda=\Vert(m,n)\Vert$, that is, $(m/\lambda,n/\lambda)$ is on the boundary of $K$. Plugging this point into 
$f(x)=1$ we get
$$\lambda=\Vert(m,n)\Vert=\sqrt[4]{4m^4-22m^3n+49m^2n^2-52mn^3+22n^4}.
$$ 
Then $\Vert (1,0)\Vert=\sqrt{2}$ and $\Vert (0,1)\Vert=\sqrt[4]{22}$. According to the algorithm, we now need to compute the norms $\Vert (m,1)\Vert$ with $m=\pm 1,\pm 2 $. The smallest of them is $\Vert (1,1)\Vert=1<\Vert (1,0)\Vert=\sqrt{2}$.
We have $\Vert (1,1)\Vert<\frac{2}{3}\Vert (0,1)\Vert$ and hence  we pass to the basis $(e^1+e^2,e^1)$. This amounts to applying matrix $\begin{bmatrix}1&1\\ 1&0\\\end{bmatrix}^{-1}$ to $K$ and hence
we  pass to $g(x,y)=f(x+y,x)=x^4-4x^3y+7x^2y^2-6xy^3+4y^4$. 
On the  next (final) iteration, we compare the norms of vectors $(m,1)$ that satisfy $|m|<2\sqrt{2}$ and pass to the basis $(e^1,e^1+e^2)$, where  $\Vert e^1+e^2\Vert=\sqrt[4]{2}$. Therefore the equation of the boundary becomes 
$h(x,y)=g(x+y,y)=x^4+x^2y^2+2y^4=1.$

It is easy to check that the Hessian of $h(x,y)$ is positive-semidefinite and hence $h(x,y)$ is convex. This implies that its level set $h(x,y)\leq 1$ is also convex and we conclude that $K$ is convex. 

Combining the two transformations, we conclude that basis $(e^1+e^2,2e^1+e^2)$ is reduced with respect to $K$. We have also shown that $\lambda_1(K)=1$ and $\lambda_2(K)=\sqrt[4]{2}$. 
\end{Ex}

Consider next a 2-dimensional sublattice of $\Z^d$. We say that its basis $(h^1,h^2)$ is reduced if it satisfies parts (1) and (2) of Definition~\ref{D:red2D}. The work that we have done above generalizes to this set up with the  only exception that we will need to use $d$ operations in order to compute a vector of the form $mh^1+h^2$, so the bound in Theorem~\ref{T:main2D} will be multiplied by $d$.


\section{Lattice Size and Basis Reduction in Dimension 3}
 
In this section we provide an algorithm for generalized basis reduction in dimension 3, thus extending the results of~\cite{KaibSch} to dimension 3. We also explain how to recover successive minima and lattice size from the obtained reduced basis of the integer lattice.

Let $K\subset\R^3$ be a convex origin-symmetric body and let $\Vert \cdot\Vert $ be the norm associated with $K$.

\begin{Def}\label{D:red3D} We say that a basis $(h^1, h^2, h^3)$ of $\Z^3$ is reduced (with respect to $K$) if
\begin{enumerate}
\item $\Vert h^1\Vert\leq  \Vert h^2\Vert \leq \Vert h^3\Vert $;
\item $\Vert h^1\pm h^2\Vert \geq  \Vert h^2\Vert$;
\item $\Vert mh^1+nh^2+h^3\Vert \geq \Vert h^3\Vert $ for all $m,n\in\Z$.
\end{enumerate}
\end{Def}

\begin{Ex} Let $P=\conv\{(0,0,0), (1,0,0), (0,1,0), (0,0,1)\}$. Let's check that the standard basis $(e^1,e^2, e^3)$ is reduced with respect to $K=\left(P+(-P)\right)^\ast$.
Indeed, we have $\Vert e^1\Vert=\Vert e^2\Vert=\Vert e^3\Vert=\Vert e^1+ e^2\Vert=1$ and $\Vert e^1-e^2\Vert=2$. Also, for $m,n\in\Z$ 
$$ \Vert mh^1+nh^2+h^3\Vert =\w_{(m,n,1)}(P)=\max\{m, n, 1\}-\min\{m, n, 0\}\geq 1=\Vert e^3\Vert,
$$
and the conclusion follows.
\end{Ex}

Note that  Definition~\ref{D:red3D} is a natural generalization of Definition~\ref{D:red2D}. As we explained in the note after that definition, the condition that $\Vert h^1\pm h^2\Vert \geq  \Vert h^2\Vert$ can be replaced with a seemingly stronger condition
that $\Vert m h^1+ h^2\Vert \geq \Vert h^2\Vert$ for all $m\in\Z$, which is of the same form as condition (3) in this definition. The same definition of a reduced basis in the case of the Euclidean norm appears in~\cite{Semaev}.

We now explain how to find the successive minima of $K$ and the lattice size of $P$ if the basis $(h^1, h^2, h^3)$ of $\Z^3$ is reduced. The situation is more complicated than in dimension two: the norms of vectors in a reduced basis do not 
necessarily match the successive minima of $K$, but we can still recover the successive minima.

\begin{Th}\label{T:main}  Let $\Vert\cdot\Vert$ be the norm associated with a convex origin-symmetric body $K$.
Let $(h^1, h^2, h^3)$ be a reduced basis of $\Z^3$ with respect to $K$ and let $\Vert\cdot\Vert$ be the norm defined by $K$.
Then $$\Vert ah^1+bh^2+ch^3\Vert \geq \Vert h^3\Vert $$ for all $(a,b,c)\in\Z^3$ with $c>0$, except, possibly, for one direction $ah^1+bh^2+ch^3\in E$, where
$$E=\{ah^1+bh^2+ch^3\mid |a| = |b| = 1\ {\rm and}\  c = 2 \}.
$$ 
This implies that $\mu_3(K)=\Vert h^3\Vert $. 

Let  $u\in E$ be the vector of the smallest norm in $E$.
Then the three successive minima $\lambda_1, \lambda_2,\lambda_3$ of $K$ are the three smallest numbers out of  $\Vert h^1\Vert $, $\Vert h^2\Vert $, $\Vert h^3\Vert$, and $\Vert u\Vert $. 
\end{Th}

\begin{Cor}\label{C:ls_computed} Let $P\subset\R^3$ be a convex body and let  $(h^1, h^2, h^3)$ be a reduced basis   of $\Z^3$ with respect to $K=(P+(-P))^\ast$. Then $\lss(P)=\mu_3(K)=\Vert h^3\Vert$.
Also, $\operatorname{ls}_1(P)$ and $\operatorname{ls}_2(P)$ are the two smallest numbers out of  $\Vert h^1\Vert $, $\Vert h^2\Vert $,  and $\Vert u\Vert $. 
\end{Cor}

In Example~\ref{E:Reeves} below we apply Theorem~\ref{T:main} and Corollary~\ref{C:ls_computed} to compute the lattice size of the Reeve tetrahedron.

\begin{proof} Fix a primitive direction $(a,b,c)\in\Z^3$ with nonzero $c$. Switching signs of the basis vectors, we can assume that $a, b, c\geq 0$.
Let $mc$ be a multiple of $c$ nearest to $a$ and $nc$ be a multiple of $c$ nearest to $b$. Then
\begin{eqnarray*}
\Vert ah^1 + bh^2 + ch^3\Vert  +  |mc-a|\Vert h^1\Vert + |nc-b|\Vert h^2 \Vert\geq  c\Vert mh^1 + nh^2 + h^3\Vert\geq c\Vert h^3\Vert,
\end{eqnarray*}
which implies
\begin{equation} \label{reduction}
	\Vert ah^1 + bh^2 + ch^3\Vert  \geq \left(c - |mc - a| - |nc - b|\right)\Vert h^3\Vert  \text{.} 
\end{equation}
Suppose first that $c$ is odd. Then
\[
	|mc - a| \le \frac{c-1}{2} \text{ and } |nc - b| \le \frac{c-1}{2} \text{.}
\]
With (\ref{reduction}), this gives 
\[
	\Vert ah^1 + bh^2 + ch^3\Vert  \ge \left(c - \frac{c - 1}{2}  - \frac{c - 1}{2}\right) \Vert h^3\Vert  = \Vert h^3\Vert \text{.}
\]
For even $c$ we have $|mc - a| \le \frac{c}{2} \text{ and } |nc - b| \le \frac{c}{2} \text{.}$
Notice that if $c\neq 2$, we cannot have equality in both of these inequalities since if this were the case $c/2$ would  divide $\gcd(a, b, c)$. 
Thus we have  $|mc-a|+|nc-b| \leq \frac{c}{2}+\frac{c-2}{2}=c-1$ and, using this in (\ref{reduction}), we get
\[
	\Vert ah^1 + bh^2 + ch^3\Vert  \ge \left(c  - c+1\right)\Vert h^3\Vert  = \Vert h^3\Vert \text{.}
\]
This leaves only the case when $c = 2$. If either $a$ or $b$ is even then one of the differences $|mc-a|$, $|nc-b|$ is zero and the other one is either zero or one and it follows that $\Vert ah^1 + bh^2 + ch^3\Vert \ge \Vert h^3\Vert $.

For the remainder of the proof, suppose that both $a$ and $b$ are odd and  $c=2$. If $\Vert h^1 + h^2\Vert  \leq \Vert h^3\Vert $, then
\begin{eqnarray*}
	\Vert ah^1 + bh^2 + 2h^3\Vert  + \Vert h^3\Vert  &\geq&\Vert ah^1 + bh^2 + 2h^3\Vert  + \Vert h^1+h^2\Vert  \\
		&\geq& \Vert (a+1)h^1 + (b+1)h^2 + 2h^3\Vert  \\
		&\geq &2 \left\Vert \frac{a+1}{2} h^1 + \frac{b+1}{2} h^2 + h^3 \right\Vert \\
		& \geq& 2 \Vert h^3\Vert,
\end{eqnarray*}
where the final inequality holds since we assumed $\Vert mh^1 + nh^2 + h^3\Vert \geq \Vert h^3\Vert $ for all $m, n \in \Z$. It follows that $\Vert ah^1 + bh^2 + 2h^3\Vert  \ge \Vert h^3\Vert $. Similarly, we get the same result  if $\Vert h^1 - h^2\Vert  \leq\Vert h^3\Vert $.

We can now assume that $\Vert h^1 \pm  h^2\Vert \geq \Vert h^3\|$. Suppose that next that $a, b\geq 3$.  Without lost of generality assume that $a\geq b$. Then
$$	\Vert ah^1 + bh^2 + 2h^3\Vert  + (a-b)\Vert   h^2\Vert + 2\Vert h^3\Vert  \ge \Vert ah^1 + ah^2\Vert  \ge a \Vert h^3\Vert,
$$
and we get 
$$\Vert ah^1 + bh^2 + 2h^3\Vert  \ge a \Vert h^3\Vert  - (a - b+ 2)\Vert h^3\Vert = (b - 2)\Vert h^3\Vert \geq \Vert h^3\Vert  \text{.}
$$
We are left with the case when $c=2$ and  $b=1$.  Using $a\geq 3$ we get
$$\Vert ah^1 + h^2 + 2h^3\Vert + \Vert h^2\Vert + (a - 2)\Vert h^3\Vert  \geq  a \Vert  h^1 + h^3\Vert  \ge a\Vert h^3\Vert .
$$
We conclude that $\Vert ah^1 + h^2 + 2h^3\Vert \geq  \Vert h^3\Vert $ and we have checked that \\ $\Vert ah^1+bh^2+ch^3\Vert \geq \Vert h^3\Vert $ for all nonzero $c$ except for the case when $|a|=|b|=1$ and $|c|=2$.
Note that if we additionally have that $\Vert h^1+h^2\Vert \leq \Vert h^3\Vert $ or $\Vert h^1-h^2\Vert \leq \Vert h^3\Vert $, then $\Vert ah^1+bh^2+ch^3\Vert \geq \Vert h^3\Vert $ holds for all nonzero $c$.

We can conclude that $\mu_3(K)=\Vert h^3\Vert $ since all the vectors whose norm is less than $\Vert h^3\Vert $ have last coordinate $0$ or $2$, and such vectors cannot form a basis of $\Z^3$.

Let $E=\{ah^1+bh^2+ch^3\mid |a| = |b| = 1\ {\rm and}\  c = 2 \}.$ 
Suppose that $u,v$ are two distinct vectors in  $E$. Then their half-sum $(u+v)/2$ has integer coordinates, is not in $E$, and has last coordinate equal to 2, so
$\Vert (u+v)/2\Vert \geq \Vert h^3\Vert$.
Hence we have
$$\Vert u\Vert +\Vert v\Vert \geq \Vert u+ v\Vert =2\Vert (u+ v)/2\Vert \geq 2 \Vert h^3\Vert,
$$
which implies that at most one vector in $E$ may have norm less than $\Vert h^3\Vert $. Define  $u\in E$ to be the vector of the smallest norm in $E$.
Any three of the four vectors in the set $\{h^1, h^2, h^3, u\}$ are linearly independent and hence we have shown that the three successive minima $\lambda_1, \lambda_2,\lambda_3$ of $K$ are the three smallest numbers on the list 
$$\{\Vert h^1\Vert , \Vert h^2\Vert , \Vert h^3\Vert , \Vert u\Vert \}.$$
 \end{proof}
 
 \begin{Rem} We would like to point out that the trick used in the proof of Theorem~\ref{T:main} (and in many other arguments in this paper) where we approximate a primitive $(a,b,c)\in\Z^3$ with $(mc,nc,c)$ does not work in dimension higher than 3.
 \end{Rem}
 
\begin{Def}\label{D:Minkowski_reduced} Let $\Vert\cdot\Vert$ be the norm that corresponds to $K$. We say that a basis $(h^1, \dots, h^d)$ of $\Z^d$ is {\it Minkowski reduced} (with respect to $K$) if $\Vert h^k\Vert=\mu_k$  for $k=1,\dots, d$.
\end{Def}

A Minkowski reduced basis is reduced and Theorem~\ref{T:main} describes how to obtain a Minkowski reduced basis of $\Z^3$ starting with a reduced one.

\begin{Def} Let $(h^1, h^2, h^3)$  be a reduced basis and let $u$ be the direction of smallest norm in $E=\{ah^1+bh^2+ch^3\mid |a| = |b| = 1\ {\rm and}\  c= 2 \}.$
Let $g^1, g^2$ be two vectors of the smaller norm in the set  $\{h^1, h^2, u\}$, written in the order of increasing norm,   and let $g^3=h^3$.
We say that $(g^1, g^2, g^3)$ is a {\it Minkowski reduced basis that corresponds to the reduced basis} $(h^1, h^2, h^3)$.
\end{Def}

\begin{Prop}\label{P:Acomputes} Let $K=(P+(-P))^\ast$ and let $(h^1,\dots, h^d)$ be a Minkowski reduced basis (with respect to $K$). Let $A$ be the matrix whose rows are the vectors of the basis $h^1, \dots, h^d$.
Then $A$ computes $\operatorname{ls}_k(P)$ for  all $k=1,\dots, d$ and, in particular, computes $\lss(P)$. 
\end{Prop}
\begin{proof}
This follows from Proposition~\ref{T:translate}, where we showed that $\operatorname{ls}_k(P)=\mu_k$.
To show that $A$ computes $\operatorname{ls}_k(P)$ we need to check that $A^{T}(\Diamond_k)\subset \mu_k K$, which is equivalent to 
checking that $\{h_1,\dots,h_k\}\subset \mu_kK$, and this holds true by the definition of the Minkowski reduced basis. 
\end{proof}

Our next goal is to develop an algorithm for finding a basis of $\Z^3$ which is reduced with respect to a given  origin-symmetric convex body $K\subset\R^3$.  
Let $(h^1,h^2, h^3)$ be a basis of $\Z^3$ such that  $\Vert h^1\Vert\leq\Vert h^2\Vert\leq\Vert h^3\Vert$ and $(h^1,h^2)$ is reduced. Let
$\min\Vert mh^1+nh^2+h^3\Vert$ over ${m,n\in\Z}$ be attained at $f=mh^1+nh^2+h^3$. 

\begin{Prop}\label{P:reduced} Let $(h^1,h^2, h^3)$ be a basis of $\Z^3$ such that  $\Vert h^1\Vert\leq\Vert h^2\Vert\leq\Vert h^3\Vert$ and $(h^1,h^2)$ is reduced. Suppose that $\Vert f\Vert\geq \Vert h^2\Vert$.
Then the basis $(h^1, h^2, f)$ is reduced.
 \end{Prop}

\begin{proof}  For integer $a,b$ we have
 $$\Vert ah^1+bh^2+ f\Vert = \Vert (a+m)h^1+(b+n)h^2+h^3 \Vert\geq \Vert f\Vert.
$$
 \end{proof}
 
 \begin{Ex}\label{E:Reeves} The Reeve tetrahedron $T_q$ is the convex hall of the columns in the matrix $\begin{bmatrix}1&0&0&1\\0&1&0&1\\0&0&1&q\end{bmatrix}$, where $q\in\mathbb{N}$.
In this example we construct a basis which is reduced with respect $K=(T_q+(-T_q))^{\ast}$ under the assumption that $q$ is odd. We have $\Vert e^1\Vert=\Vert e^2\Vert=1$ and $\Vert e^3\Vert=q$. Since $\Vert e^1\pm e^2\Vert=2$, we conclude that $(e^1,e^2)$ is reduced.
Denote $h=-\frac{q-1}{2}e^1-\frac{q-1}{2}e^2+e^3$. Then
$$\Vert h\Vert=\max\left\{-\frac{q-1}{2},1\right\}-\min\left\{-\frac{q-1}{2},1\right\}=\frac{q+1}{2}.
$$
For $m,n\in\Z$ we have
$$\Vert{mh^1+nh^2+h^3}\Vert=\max\{m,n,1,m+n+q\}-\min\{m,n,1,m+n+q\}\geq 1-m.
$$
We also have $\Vert{mh^1+nh^2+h^3}\Vert\geq m+n+q-n=m+q$, and taking the half-sum of the two inequalities we conclude that $\Vert{mh^1+nh^2+h^3}\Vert\geq \frac{q+1}{2}=\Vert h\Vert$.
By Proposition~\ref{P:reduced}  basis $(e^1,e^2, h)$ is reduced. It is easy to check that $$\Vert \pm e^1\pm e^2+2e^3\Vert=q+2,$$ and hence by Theorem~\ref{T:main} and Corollary~\ref{C:ls_computed} we conclude that 
$\lambda_1(K)=\lambda_2(K)=1$ and $\lambda_3(K)=\mu_3(K)=\lss(T_q)=\frac{q+1}{2}$.
\end{Ex}

Proposition~\ref{P:reduced} implies that the following algorithm terminates in a reduced basis.
 
\begin{Alg}\label{A:alg3}  Let $\Vert\cdot\Vert$ be the norm associated with a convex origin-symmetric body $K$.
Order the basis vectors $(h^1,h^2, h^3)$ of $\Z^3$ so that $\Vert h^1\Vert\leq\Vert h^2\Vert\leq\Vert h^3\Vert$ and reduce $(h^1,h^2)$.
At each step, find $f$ such that the minimum of $\Vert mh^1+nh^2+h^3\Vert$ over $m, n\in\Z$ is attained at~$f=mh^1+nh^2+h^3$.
If $\Vert f\Vert< \Vert h^2\Vert$, pass to the basis $(h^1,h^2,f)$ and repeat the step. If $\Vert f\Vert\geq \Vert h^2\Vert$  output a reduced basis  
$(h^1,h^2,f)$.
\end{Alg}

This algorithm terminates after finitely many iterations (see the argument after Algorithm~\ref{A:alg1}).  In order to be able to estimate the number of iterations, we will modify the algorithm, similarly to how we did this in dimension 2.  In preparation for this we prove the next two propositions.

\begin{Prop}\label{P:halfsum} Let  $(h^1, h^2)$  be a reduced basis and let $\Vert h^1\Vert\leq \Vert h^2\Vert\leq \Vert h^3\Vert$. Suppose that $\frac{\Vert h^1\Vert+\Vert h^2\Vert}{2}\leq \Vert f\Vert< \Vert h^2\Vert$. Then Algorithm~\ref{A:alg3} will terminate in at most two additional iterations.
\end{Prop}
\begin{proof}
At the initial step, we pass from $(h^1,h^2,h^3)$ to $(h^1, f, h^2)$.
 At the next step, let $$\min\limits_{a,b\in\Z}\Vert ah^1+b(mh^1+nh^2+h^3)+h^2\Vert$$ be attained at $g=(a+bm)h^1+(bn+1)h^2+bh^3$. 
If  $b=0$ then $\Vert g\Vert\geq \Vert h^2\Vert$; if $|b|=1$ then $\Vert g\Vert\geq \Vert f\Vert$, so the algorithm terminates at this step. 

Next suppose that $|b|\geq 2$ and 
let  $bk$ be the multiple of $b$ closest to $a$. Then 
\begin{equation*}
\Vert g\Vert+|a-bk|\Vert h^1\Vert+\Vert h^2\Vert\geq |b|\Vert (m+k)h^1+nh^2+h^3\Vert\geq |b|\Vert f\Vert, 
\end{equation*}
which for $|b|\geq 4$ implies 
\begin{eqnarray*}
\Vert g\Vert&\geq &|b|\Vert f\Vert-|a-bk|\Vert h^1\Vert-\Vert h^2\Vert\geq |b|\Vert f\Vert -\frac{|b|}{2}\Vert h^1\Vert- \Vert h^2\Vert\\
&\geq& |b|\frac{\Vert h^1\Vert+\Vert h^2\Vert}{2}-\frac{|b|}{2}\Vert h^1\Vert- \Vert h^2\Vert=\frac{|b|-2}{2}\Vert h^2\Vert\geq \Vert h^2\Vert.
\end{eqnarray*}
If $|b|=3$ then  $|a-bk|\leq 1$ and we get 
\begin{eqnarray*}
\Vert g\Vert&\geq &|b|\Vert f\Vert-|a-bk|\Vert h^1\Vert-\Vert h^2\Vert\geq 3\Vert f\Vert -\Vert h^1\Vert- \Vert h^2\Vert\geq \Vert f\Vert.
\end{eqnarray*}

Let next $|b|=2$. If the algorithm has not terminated at the previous step, we have $\Vert g\Vert<\Vert f\Vert$, so we 
are now minimizing over $A,B\in\Z$ the norm of a vector of the form $Ah^1+Bg+f$. If we rewrite this vector in terms of $h^1, h^2,$ and $h^3$ then the coefficient 
of $h^3$ is $Bb+1$, which is odd since  $|b|=2$. Consider $p,q,r\in \Z$ where $r$ is odd and positive.  Then, working with multiples of $r$ closest to $p$ and $q$,  we get
$$\Vert ph^1+qh^2+rh^3\Vert\geq r\Vert f\Vert-\frac{r-1}{2}\Vert h^1\Vert-\frac{r-1}{2}\Vert  h^2\Vert\geq \Vert f\Vert,$$
and hence the algorithm terminates at this step.
%
%
\end{proof}

Let $\min\{\Vert h^1\pm h^2\pm 2h^3\Vert\}$  be attained at~$g$.
The following statement is the three-dimensional version of Proposition~\ref{P:exceptional2D}.

\begin{Prop}\label{P:exceptionalS}  Let  $(h^1, h^2)$  be reduced and $\Vert h^1\Vert\leq \Vert h^2\Vert\leq \Vert h^3\Vert$. 
Suppose that $\Vert f\Vert<\frac{\Vert h^1\Vert+\Vert h^2\Vert}{2}$. Let $\Vert f\Vert>\frac{19}{20}\Vert h^3\Vert$ and  $\Vert g\Vert>\frac{9}{10}\Vert h^3\Vert$.
Then $$\Vert ah^1+bh^2+ch^3\Vert\geq \Vert h^3\Vert$$ for all primitive $(a,b,c)\in\Z^3$ except for the ones for which $(|a|,|b|,|c|)$ or $(|b|,|a|,|c|)$ is on the list below:
\begin{gather*}
(1,1,2),(2,1,3),(2,2,3),(3,2,4),(3,2,5),(1,0,0), (1,0,1), (1,1,0), \\
(1,1,1), (2,0,1), (2,1,0), (2,1,1), (2,1,2), (2,2,1),(3,1,1), (3,1,2), \\ 
(3,2,1), (3,2,2), (4,1,2), (4,2,1), (4,2,3), (4,3,2), (5,2,2), (5,2,3), (5,3,2).
\end{gather*}
We denote this set of exceptional directions by $S$. Note that $S$ has 160 directions if we do not distinguish between directions $(a,b,c)$ and $(-a,-b,-c)$. 
\end{Prop}

\begin{proof}
We have $\Vert h^2\Vert>\Vert f\Vert>\frac{19}{20}\Vert h^3\Vert$ and 
$$\Vert h^1\Vert>2\Vert f\Vert- \Vert h^2\Vert\geq \frac{38}{20}\Vert h^3\Vert-\Vert h^3\Vert = \frac{9}{10}\Vert h^3\Vert.
$$
Let $(a,b,c)\in\Z^3$ be primitive. Without loss of generality we can assume that $a,b,c\geq 0$. In what follows we will use $\Vert h^1\Vert, \Vert h^2\Vert\geq \frac{9}{10}\Vert h^3\Vert$, although we have a stronger inequality on $\Vert h^2\Vert$, to make the argument symmetric in $a$ and $b$. For the same reason, we will not be using the inequality $\Vert h^1\Vert\leq \Vert h^2\Vert$, so without loss of generality we can assume $a\geq b$.

Suppose first that $c> a\geq b$.  We then have
\begin{eqnarray*}
&&\Vert ah^1+bh^2+ch^3\Vert+a\Vert h^1\Vert+b\Vert h^2\Vert\geq c\Vert h^3\Vert;\\
&&\Vert ah^1+bh^2+ch^3\Vert+(c-a)\Vert h^1\Vert+b\Vert h^2\Vert\geq c\Vert h^1+h^3\Vert\geq c\Vert f\Vert> \frac{19}{20}c\Vert h^3\Vert;\\
&& \Vert ah^1+bh^2+ch^3\Vert+ (c-a)\Vert h^1\Vert+(c-b)\Vert h^2\Vert \geq c\Vert f\Vert\geq \frac{19}{20} c \Vert h^3\Vert.
\end{eqnarray*}
Working with each of the above inequalities separately we get
\begin{eqnarray*}
&&\Vert ah^1+bh^2+ch^3\Vert\geq c\Vert h^3\Vert-a\Vert h^1\Vert-b\Vert h^2\Vert\geq (c-a-b)\Vert h^3\Vert;\\
&&\Vert ah^1+bh^2+ch^3\Vert\geq \left(\frac{19}{20}c-(c-a)-b\right)\Vert h^3\Vert= \left(a-b-\frac{c}{20}\right)\Vert h^3\Vert;\\
&& \Vert ah^1+bh^2+ch^3\Vert\geq  \left(\frac{19}{20} c - (c-a)-(c-b)\right) \Vert h^3\Vert= \left(a+b-\frac{21}{20} c \right) \Vert h^3\Vert.
\end{eqnarray*}
Hence we can conclude that $\Vert ah^1+bh^2+ch^3\Vert\geq   \Vert h^3\Vert$ unless
\begin{equation}\label{e:system_abc1}
c\leq a+b,\ \  \ a<b+1+\frac{c}{20},\ \ \ a+b<\frac{21}{20}c+1.
\end{equation}
 Combining these inequalities together with $c\geq a\geq b$ we get 
 $$2b-1-\frac{c}{20}\leq \frac{20}{21}(2b-1)\leq \frac{20}{21}(a+b-1)\leq c \leq a+b\leq 2b+1+\frac{c}{20},
 $$
 where  the first inequality follows from second and third in this chain.

We conclude that  $|c-2b|\leq 1+\frac{c}{20}$. We also have $0\leq a-b\leq 1+\frac{c}{20}$. Hence 
\begin{gather*}
\Vert ah^1+bh^2+ch^3\Vert+(a-b)\Vert h^1\Vert+|c-2b|\Vert h^3\Vert\geq b\Vert g\Vert\geq \frac{9}{10}b\Vert h^3\Vert, \ \ {\rm so}\\
\Vert ah^1+bh^2+ch^3\Vert\geq \left(\frac{9}{10}b-(a-b)-|c-2b|\right)\Vert h^3\Vert\geq\left(\frac{9}{10}b-2-\frac{c}{10}\right)\Vert h^3\Vert,
\end{gather*}
and  we can assume that $\frac{9}{10}b-2-\frac{c}{10}<1$, that is, $9b-29\leq c$. From $c\leq 2b+1+\frac{c}{20}$ we get $c\leq\frac{20}{19}(2b+1)$ and hence
$9b-29\leq \frac{20}{19}(2b+1)$, 
which implies that $b\leq 4$. It follows that $c\leq 9$ and $a=b$ or $b+1$.
Note that we can also assume that 
\begin{equation}\label{e:abs_sum}
\frac{9}{10}b-(a-b)-|c-2b|<1.
\end{equation}
We use Maple~\cite{Maple} to find all  primitive triples  $(a,b,c)\in\Z^3$ where $0\leq b\leq 4$, $1\leq c\leq 9$, $a=b$ or $b+1$  that satisfy  $c> a$ together with  (\ref{e:system_abc1}) and (\ref{e:abs_sum}).
The output is $\{ (1,1,2),  (2,1,3), (2,2,3), (3,2,4), (3,2,5), (4,3,7)\}$, and we can eliminate the last vector on this list using
$$\Vert 4h^1+3h^2+7h^3\Vert\geq 4\Vert g\Vert -\Vert h^2\Vert- \Vert h^3\Vert\geq \left(4\cdot\frac{9}{10}-2\right)\Vert h^3\Vert> \Vert h^3\Vert.
$$

We next take care of  the  remaining case  $a\geq c,b$. We have 
\begin{eqnarray*}
&&\Vert ah^1+bh^2+ch^3\Vert+b\Vert h^2\Vert+c\Vert h^3\Vert\geq a\Vert h^1\Vert> \frac{9}{10}a\Vert h^3\Vert;\\
&& \Vert ah^1+bh^2+ch^3\Vert+ b\Vert h^2\Vert+(a-c)\Vert h^3\Vert \geq a\Vert f\Vert> \frac{19}{20} a\Vert h^3\Vert;\\
&&\Vert ah^1+bh^2+ch^3\Vert+(a-b)\Vert h^2\Vert+c\Vert h^3\Vert\geq a\Vert h^1+h^2\Vert\geq a\Vert h^2\Vert> \frac{9}{10}a\Vert h^3\Vert;\\
&& \Vert ah^1+bh^2+ch^3\Vert+ (a-b)\Vert h^2\Vert+(a-c)\Vert h^3\Vert \geq a\Vert f\Vert> \frac{19}{20} a \Vert h^3\Vert.
\end{eqnarray*}

Hence we can conclude  $\Vert ah^1+bh^2+ch^3\Vert\geq \Vert h^3\Vert$ unless $a,b,$ and $c$ satisfy the following system of inequalities
\begin{equation}\label{e:system_abc}
\frac{9}{10}a<b+c+1,\ \ c<b+1+\frac{a}{20},\ \ b<c+1+\frac{a}{10},\ \ b+c<1+\frac{21}{20}a,
\end{equation}
which we assume  to hold true from now on. 

The second and third inequality together imply that $|b-c|< 1+\frac{a}{10}$.
We also have $\frac{9}{10}a-c-1<b<c+\frac{a}{10}+1$, which implies that $a<\frac{5}{2}(c+1)$
and hence $|b-c|< 1+\frac{c+1}{4}$.

Let $kc$ be the multiple of $c$ closest to $a$. We get
\begin{equation*}
\Vert ah^1+bh^2+ch^3\Vert+|a-kc|\Vert h^1\Vert+|b-c|\Vert h^2\Vert\geq \frac{19}{20}c\Vert h^3\Vert,
\end{equation*}
and hence
\begin{equation*}
\Vert ah^1+bh^2+ch^3\Vert \geq \left(\frac{19}{20}c -\frac{c}{2}-\frac{c+1}{4}-1\right)\Vert h^3\Vert\geq \Vert h^3\Vert,
\end{equation*}
unless $c\leq 11$, which we  now assume. This implies that $a\leq 29$ and $b\leq 14$.

Approximating $(a,b,c)$ first with $(2b,b,b)$ and next with $(2c,c,c)$ we get
\begin{eqnarray*}
\Vert ah^1+bh^2+ch^3\Vert &\geq &\left(\frac{19}{20}b -|a-2b|-|b-c|\right)\Vert h^3\Vert,\\
\Vert ah^1+bh^2+ch^3\Vert &\geq &\left(\frac{19}{20}c -|a-2c|-|b-c|\right)\Vert h^3\Vert,\\
\end{eqnarray*}
so  we can conclude  that $\Vert ah^1+bh^2+ch^3\Vert\geq \Vert h^3\Vert$ unless
\begin{equation}\label{e:abs_val} 
\frac{19}{20}b -|a-2b|-|b-c|<1 \ \ {\rm and}\  \  \frac{19}{20}c -|a-2c|-|b-c|<1.
\end{equation}
We next use Maple to find all primitive  triples $(a,b,c)\in\Z^3$ where $1\leq a\leq 29$, $0\leq b\leq 14$,  and $0\leq c\leq 11$, that satisfy $a\geq b$ and $a\geq c$ together with (\ref{e:system_abc}) and (\ref{e:abs_val}).
The output is 
\begin{align*}
(1,0,0), (1,0,1), (1,1,0), (1,1,1), (2,0,1), (2,1,0), (2,1,1), (2,1,2), (2,2,1),(3,1,1),\\ 
(3,1,2), (3,2,1), (3,2,2), (4,1,2), (4,2,1), (4,2,3), (4,3,2), (5,2,2), (5,2,3), (5,3,2).
\end{align*}
\end{proof}

The algorithm below starts with a basis $(h^1,h^2,h^3)$ of $\Z^3$ and outputs a basis which is reduced with respect to a given origin-symmetric body $K$. 
By Theorem~\ref{T:main} one can then easily recover the successive minima of $K$ comparing the norms of the output vectors and the norm of the vector $u$, defined in that theorem. 
If we start with a convex body $P\subset\R^3$ and apply the algorithm to $K=\left(P+(-P)\right)^{\ast}$ the algorithm computes the lattice size of $P$.

\begin{Alg}\label{A:alg4}  Let $\Vert\cdot\Vert$ be the norm associated with a convex origin-symmetric body $K$.
Order the basis vectors $(h^1,h^2, h^3)$ of $\Z^3$ so that $\Vert h^1\Vert\leq\Vert h^2\Vert\leq\Vert h^3\Vert$ and reduce $(h^1,h^2)$ using  Algorithm~\ref{A:alg2}.
At each step, find $f$ such that the minimum of $\Vert mh^1+nh^2+h^3\Vert$ over $m, n\in\Z$ is attained at~$f=mh^1+nh^2+h^3$.
\begin{itemize}
\item[(1)] If $\Vert f\Vert\geq \Vert h^2\Vert$ output $(h^1, h^2, f)$, which is reduced by Proposition~\ref{P:reduced}.
\item[(2)] If $ \frac{1}{2}\left(\Vert h^1\Vert+\Vert h^2\Vert\right)\leq \Vert f\Vert<\Vert h^2\Vert$, pass to $(h^1, f, h^2)$ and repeat the step.
By Proposition~\ref{P:halfsum} we will  get a reduced basis after at most two additional iterations.
\item[(3)] If $\Vert f\Vert\leq \frac{19}{20}\Vert h^3\Vert$  and $\Vert f\Vert< \frac{1}{2}\left(\Vert h^1\Vert+\Vert h^2\Vert\right)$ repeat the step with the basis $(h^1,f,  h^2)$. \vspace{.1cm}
\item[(4)] If $\frac{19}{20}\Vert h^3\Vert<\Vert f\Vert< \frac{1}{2}\left(\Vert h^1\Vert+\Vert h^2\Vert\right)$  let $\min\{\Vert h^1\pm h^2\pm 2h^3\Vert\}$  be attained at~$g$.
	\begin{itemize}
		\item[(i)] If $\Vert g\Vert\leq \frac{9}{10}\Vert h^3\Vert$ then pass to $(h^1, g, h^3)$ and repeat the step.
	\vspace{.1cm}
	\item[(ii)] If  $\Vert g\Vert>\frac{9}{10}\Vert h^3\Vert$ then the exceptional set $S$ from Proposition~\ref{P:exceptionalS} contains a Minkowski reduced basis. Let $h^1$ be a vector of smallest norm in $S$, let $h^2\in S$ be a vector of smallest norm linear independent of $h^1$ and, finally, let  $h^3\in S$ be a vector of smallest norm such that $(h^1,h^2, h^3)$ is a basis of $\Z^3$. Output a Minkowski reduced basis $(h^1,h^2,h^3)$.
	\end{itemize}
\end{itemize}
\end{Alg}

\begin{Th}
The number of iterations in Algorithm~\ref{A:alg4} is at most $$\log_{\frac{60}{59}}\frac{\Vert h^1\Vert+\Vert h^2\Vert+\Vert h^3\Vert}{\mu_1+\mu_2+\mu_3}+2.$$
\end{Th}
\begin{proof}
We track how $\Vert h^1\Vert+\Vert h^2\Vert+\Vert h^3\Vert$ changes throughout the iterations.
In (3) we have $\Vert f\Vert\leq \frac{19}{20}\Vert h^3\Vert$  and we pass from $(h^1,h^2,h^3)$ to $(h^1,h^2,f)$.
The sum of norms goes down by at least
$$\Vert h^3\Vert-\Vert f\Vert\geq \Vert h^3\Vert-\frac{19}{20}\Vert h^3\Vert=\frac{\Vert h^3\Vert}{20}\geq \frac{1}{60} \left(\Vert h^1\Vert+\Vert h^2\Vert+\Vert h^3\Vert\right).
$$

In 4(i) we have 
$$\Vert h^2\Vert - \Vert g\Vert\geq \frac{19}{20}\Vert h^3\Vert-\frac{9}{10}\Vert h^3\Vert=\frac{\Vert h^3\Vert}{20}\geq \frac{1}{60} \left(\Vert h^1\Vert+\Vert h^2\Vert+\Vert h^3\Vert\right).
$$
 We also have at most two iterations in (2), and the bound follows.
\end{proof}

Our next goal is to estimate how long each iteration in Algorithm~\ref{A:alg4} takes, that is, we need to discuss how to minimize $\Vert mh^1+nh^2+h^3\Vert$ over $m,n\in\Z$.

\begin{Prop}\label{P:mn} Let  $\Vert h^1\Vert\leq \Vert h^2\Vert\leq \Vert h^3\Vert$  and let $(h^1,h^2)$ be reduced.
Then for $(m,n)\in\Z^2$ we have  $\Vert mh^1+nh^2+h^3\Vert \geq \Vert h^3\Vert $ unless
\begin{enumerate}
\item $|m|\geq |n|$, $|n| < \frac{2\Vert h^3\Vert }{\Vert h^2\Vert }$, and $|m|< \frac{2\Vert h^3\Vert + |n|\Vert h^2\Vert }{\Vert h^1\Vert }$, or
\item $|n|>|m|$, $|m|< \frac{2\Vert h^3\Vert }{\Vert h^1\Vert }$, and $|n|< \frac{2\Vert h^3\Vert + |m|\Vert h^1\Vert }{\Vert h^2\Vert }$.

\end{enumerate}
\end{Prop}
\begin{proof} Without loss of generality we can assume that $m\geq n\geq 0$.
If $n\geq  \frac{2\Vert h^3\Vert }{\Vert h^2\Vert }$ we get
$$m\Vert h^2\Vert \leq m\Vert h^1+h^2\Vert \leq \Vert mh^1+nh^2+h^3\Vert +(m-n)\Vert h^2\Vert +\Vert h^3\Vert ,
$$
which implies
$$\Vert mh^1+nh^2+h^3\Vert \geq n\Vert h^2\Vert -\Vert h^3\Vert \geq 2\Vert h^3\Vert -\Vert h^3\Vert =\Vert h^3\Vert .
$$
If $m\geq \frac{2\Vert h^3\Vert + |n|\Vert h^2\Vert }{\Vert h^1\Vert }$ we get
$$\Vert mh^1+nh^2+h^3\Vert \geq |m|\Vert h^1\Vert -|n|\Vert h^2\Vert -\Vert h^3\Vert \geq  2\Vert h^3\Vert + |n|\Vert h^2\Vert -|n|\Vert h^2\Vert -\Vert h^3\Vert =\Vert h^3\Vert .
$$

Observe that if, say, $|m|\geq |n|$, then we have $|n| < \frac{2\Vert h^3\Vert }{\Vert h^2\Vert }$ and $|m| < \frac{4\Vert h^3\Vert }{\Vert h^1\Vert }$,  and hence
the size of the exceptional set where we may have $\Vert mh^2+nh^2+h^3\Vert < \Vert h^3\Vert $ is  $\mathcal{O}\left(\frac{\Vert h^3\Vert ^2}{\Vert h^1\Vert \Vert h^2\Vert }\right)$.
\end{proof}

It follows from this proposition that each iteration in Algorithm~\ref{A:alg4} takes at most $\mathcal{O}\left(\frac{\Vert h^3\Vert ^2}{\mu_1\mu_2}\right)$ operations.
Similarly to what we did in the proof of Theorem~\ref{T:main2D}, we show next that we need that many operations only at the final iteration of the algorithm, while at all the previous iterations the number
of directions that we need to check is smaller.

Assume first that we are not at the final iteration, that is, $\Vert f\Vert<\Vert h^2\Vert$. Then to find the minimum of $\Vert mh^1+nh^2+h^3\Vert$ we do not need to consider directions with \\ $\Vert mh^1+nh^2+h^3\Vert\geq \Vert h^2\Vert$.
Suppose first that $|m|\leq |n|$. Then 
$$\Vert mh^1+nh^2+h^3 \Vert +(|n|-|m|)\Vert h^1\Vert+\Vert h^3\Vert\geq |n|\Vert h^1\pm h^2\Vert\geq |n|\Vert h^2\Vert
$$
and hence $\Vert mh^1+nh^2+h^3 \Vert \geq |m|\Vert h^2 \Vert -\Vert h^3 \Vert\geq \Vert h^2 \Vert$ unless $|m|<1+\frac{\Vert h^3 \Vert}{\Vert h^2 \Vert}$.
Now for each such value of $m$ if we also have $|n|\geq \frac{\Vert mh^1+h^3\Vert}{\Vert h^2\Vert}+1$ then
$$\Vert mh^1+nh^2+h^3 \Vert\geq |n|\Vert h^2\Vert-\Vert mh^1+h^3\Vert \geq\Vert h^2\Vert.
$$
For $|n|\leq \frac{\Vert mh^1+h^3\Vert}{\Vert h^2\Vert}-1$ we have
$$\Vert mh^1+nh^2+h^3 \Vert\geq \Vert mh^1+h^3\Vert -|n|\Vert h^2\Vert\geq\Vert h^2\Vert.
$$
Hence the only values of $n$ for which the norm could be less than $\Vert h^2\Vert$ are $\pm\left\lfloor \frac{\Vert mh^1+h^3\Vert}{\Vert h^2\Vert}\right\rfloor$ and 
$\pm\left\lceil  \frac{\Vert mh^1+h^3\Vert}{\Vert h^2\Vert}\right\rceil$.
We have checked that the number of directions that we need to check under the assumption that $|m|\leq |n|$ is at most $c\frac{\Vert h^3\Vert}{\Vert h^2\Vert}\leq c\frac{\Vert h^3\Vert}{\Vert \mu_1\Vert}$ for some constant $c$.
Suppose next that $|m|>|n|$. Then 
$$\Vert mh^1+nh^2+h^3 \Vert +(|m|-|n|)\Vert h^2\Vert+\Vert h^3\Vert\geq |m|\Vert h^1+h^2\Vert\geq |m|\Vert h^2\Vert
$$
and hence 
$$\Vert mh^1+nh^2+h^3 \Vert\geq n\Vert h^2\Vert-\Vert h^3\Vert\geq \Vert h^2\Vert$$
unless $|n|<\frac{\Vert h^3\Vert}{\Vert h^2\Vert}$+1. Now for each such fixed $n$ we have 
$$\Vert mh^1+nh^2+h^3 \Vert\geq \Vert nh^2+h^3 \Vert -\Vert h^1 \Vert\geq \Vert h^2\Vert
$$
unless $|m|>\frac{\Vert nh^2+h^3\Vert-\Vert h^2\Vert}{\Vert h^1\Vert}$ and 
$$\Vert mh^1+nh^2+h^3 \Vert\geq |m|\Vert h^1 \Vert- \Vert nh^2+h^3\Vert\geq \Vert h^2\Vert
$$
unless $|m|<\frac{\Vert nh^2+h^3\Vert+\Vert h^2\Vert}{\Vert h^1\Vert}$. That is, the remaining values of $m$ satisfy  
$$\frac{\Vert nh^2+h^3\Vert-\Vert h^2\Vert}{\Vert h^1\Vert}<|m|<\frac{\Vert nh^2+h^3\Vert+\Vert h^2\Vert}{\Vert h^1\Vert}$$
and  we have $c\frac{\Vert h^3\Vert}{\Vert h^2\Vert}\frac{\Vert h^2\Vert}{\Vert h^1\Vert}=c\frac{\Vert h^3\Vert}{\Vert h^1\Vert}\leq c\frac{\Vert h^3\Vert}{\mu_1}$ directions to check in this case as well.
If we check all these directions and all the corresponding norms exceed $\Vert h^2\Vert$, it means that we are in the final iteration and we need to check all the directions of Proposition~\ref{P:mn}. 

As in Theorem~\ref{T:main2D}, we assume that we have access to the norm oracle. We proved the following theorem.
\begin{Th}
The number of operations in Algorithm~\ref{A:alg4} is at most $$\mathcal{O}\left(\frac{\Vert h^3\Vert}{\mu_1}\left(\log\frac{\Vert h^1\Vert+\Vert h^2\Vert+\Vert h^3\Vert}{\mu_1+\mu_2+\mu_3}\right)+\frac{\Vert h^3\Vert^2}{\mu_1\mu_2}\right).$$
\end{Th}

 \subsection*{Acknowledgments} 
We are grateful to the anonymous referee of an earlier version of this paper for explaining  to us the connection of the lattice size to the successive minima and for encouraging us to make  our arguments work for an arbitrary convex body $P\subset\R^3$.
We would also like to thank the anonymous referee of the current version for making sure that the motivation is carefully explained and that an ample amount of examples is provided.


%
%
%

\end{document}